\documentclass{iopart}
\usepackage[utf8]{inputenc}
\usepackage[T1]{fontenc}
\usepackage[greek,main=english]{babel}
\languageattribute{greek}{ancient}
\usepackage{alphabeta}
\usepackage{iopams}
\usepackage{upgreek}
\usepackage{booktabs}
\usepackage{tikz-cd}
\usepackage{amsmathe}
\usepackage{amsthm}
\usepackage[%
	backend=biber,
	style=numeric,
	hyperref=true,
	sorting=none
]{biblatex}

\newcommand{\E}{\mathrm{e}}
\newcommand{\I}{\mathrm{i}}
\newcommand{\D}{\mathrm{d}}
\newcommand{\PI}{\uppi}
\newcommand{\cord}{\operatorname{cord}}
\newcommand{\Sin}{\operatorname{Sin}}
\newcommand{\Cos}{\operatorname{Cos}}
\newcommand{\Arcsin}{\operatorname{Arcsin}}
\newcommand{\Arccos}{\operatorname{Arccos}}
\newtheorem{definition}{Definition}
\newtheorem{theorem}{Theorem}

\begin{document}

\title{The term `angle' in the international system of units}
\author{Michael P. Krystek}
\address{Physikalisch-Technische Bundesanstalt,\\ Bundesallee 100, D-38116 Braunschweig, Germany}
\eads{\mailto{Michael.Krystek@ptb.de}}
%\submitto{Metrologia}

\begin{abstract}
	The concept of an angle is one that often causes difficulties in metrology. These are partly caused by a confusing mixture of several mathematical terms, partly by real mathematical difficulties and finally by imprecise terminology. The purpose of this publication is to clarify misunderstandings and to explain why strict terminology is important. It will also be shown that most misunderstandings regarding the `radian' can be avoided if some simple rules are obeyed.
\end{abstract}


\section{Introduction}

The concept of an angle is one that often causes difficulties in metrology. These are partly caused by a confusing mix of several mathematical terms, partly by real mathematical difficulties, and finally also by inaccurate terminology.

Confusion arises {e.\,g.}\ from using the term `angle' to refer to more or less connected, but nonetheless different mathematical objects, such as a pair of intersecting straight lines, a pair of rays with a common point, a cluster of rays bounded by two rays, a circular arc, a convex part of the plane, a sector of a circle, or a rotation.

The mathematical difficulties are mainly due to the fact that statements that apply to lengths cannot be directly transferred to angular magnitudes, because the structures of angles and straight lines are different.

Inaccurate terminology causes misunderstandings. The term `angle' is used with different meanings. In mathematics, the term `plane angle', for example, refers to a geometric object (a figure, such as a straight-line segment). To geometric objects at least one magnitude can be assigned (to a straight-line segment, for example, the magnitude called `length'). The magnitude assigned to an angle is often referred to as `angle'. Moreover, the numerical value associated with this magnitude is usually also called `angle'. The usage of the same term with different meanings might be acceptable in daily life, but not in basic documents of metrology.

In the following sections, we will deal in a mathematically rigorous way with the various terms associated with the concept of an angle. The aim is to clarify misunderstandings and to show why a strict terminology is important.

\section{The concept of a plane angle}

Since ancient times, several definitions of the term `plane angle' have been given. But as we are not dealing with the history of mathematics, we go directly to the definition that is generally accepted in mathematics today (see Hilbert \cite{Hilbert1902}, for example).

\begin{quote}
	\begin{definition}
		A `plane angle', denoted by $\angle\,pq$, is a geometric figure, formed by a pair of distinct rays $p$ and $q$, called the `sides' of the angle, that originate from a common point $O$, called the `vertex' of the angle.
	\end{definition}
	\par\vspace*{0.5\baselineskip}
	\textbf{Note}\quad A plane angle may also be defined by a triple of points, then written in the form $\angle\,POQ$ instead of $\angle\,pq$, where $O$ denotes the vertex, $P$ any point on the side $p$, and $Q$ any point on the side $q$.
	\begin{figure}[htp]
		\centering
		\begin{tikzpicture}[x=1mm,y=1mm,scale=0.75]
			\draw[thick] (30:35) -- (0,0) -- (-30:35);
			\draw[right] (30:35) node {$p$};
			\draw[right] (-30:35) node {$q$};
			\filldraw (30:26) circle[radius=0.6mm];
			\draw[above] (30:26) node {$P$};
			\filldraw (0,0) circle[radius=0.6mm];
			\draw[left] (0,0) node {$O$};
			\filldraw (-30:16) circle[radius=0.6mm];
			\draw[above] (-30:16) node {$Q$};
			\draw[dashed] (0,0) circle[radius=20];
			\draw (-30:20) arc (-30:30:20);
			\draw[right] (20,0) node {$s$};
		\end{tikzpicture}
		\caption{A plane angle.}
	\label{fig1}
	\end{figure}
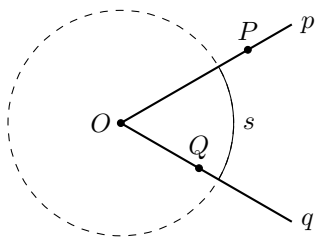
\end{quote}

The definition of a plane angle does not require an \emph{ordered} pair of rays, {i.\,e.}\ it does not distinguish between $\angle\,pq$ and $\angle\,qp$, or between $\angle\,POQ$ and $\angle\,QOP$. Following Euclid \cite{Euclid1968}, this is still common in Euclidean geometry of the plane, in which plane angles are not considered to be signed.

The sectors of a circle created by an angle are usually convex figures. This property guarantees uniqueness by associating to an angle always the shorter of the two possible arcs of a circle around its vertex, which is denoted by $s$ in Figure \ref{fig1}. The convexity of a sector excludes so-called `reflex angles'.

Due to the distinctness of the sides, the arc $s$ has always a non-zero length, {i.\,e.}\ the 'zero angle' does not exist. In contrast to the definition given by Hilbert \cite{Hilbert1902}, however, the sides of an angle can lie on the same straight line (then they, of course, have opposite orientation), {i.\,e.}\ the `straight angle' is not excluded.

\section{Angular magnitudes}

In geometry, equality is expressed by congruency. Two geometrical objects are said to be congruent if, and only if, each can be transformed into the other by an isometry. Congruency is an equivalence relation, which provides a partition of a set of geometrical objects into disjoint equivalence classes.

To every angle a magnitude can be assigned. All congruent angles, {i.\,e.}\ those which belong to the same equivalence class, have the \emph{same} magnitude. Thus we can equate the magnitude of all congruent angles with their equivalence class.

\begin{quote}
	\begin{definition}
		The equivalence class
		\begin{equation*}
			[\angle\,pq]_{\cong}:=
			\left\{\angle\,uv\,\mid\,\angle\,uv\cong\angle\,pq\right\}\,,
		\end{equation*}
		{i.\,e.}\ the set of all angles congruent to $\angle\,pq$, is called their `angular magnitude'.
	\end{definition}
	\par\vspace*{0.5\baselineskip}
	\textbf{Note}\quad The angle denoted here by $\angle\,pq$ is called a `representative' of the equivalence class $[\angle\,pq]_{\cong}$. Any element of an equivalence class can be chosen as its representative.
\end{quote}

Traditionally, the `angular magnitude' is often referred to just as `angle'. There is no harm to do this (except in standardization, where strictness is important), provided it is clear from the context whether reference is made to a geometric object or its assigned magnitude. Most misunderstandings can be avoided if a clear distinction between angles as geometric objects and their magnitudes is made.

All angular magnitudes form a set, having the cardinality of a continuum, on which the structure of an ordered commutative semigroup can be established by introducing a binary operation between the elements of the set.

\begin{quote}
	\begin{theorem}\label{Def3}
		Let $\mathcal{A}$ be a set of angular magnitudes and let $\oplus:\mathcal{A}\times\mathcal{A}\mapsto\mathcal{A}$ be a binary operation, defined by
		\begin{equation*}
			[\angle\,pr]_{\cong}\oplus[\angle\,rq]_{\cong}=[\angle\,pq]_{\cong}\,,
		\end{equation*}
		then $(\mathcal{A},\oplus)$ is an ordered, commutative and cancellative semigroup.
	\end{theorem}
	\par\vspace*{0.5\baselineskip}
	\textbf{Note 1}\quad The operation $\oplus$ is called `addition of angular magnitudes'. This addition is done modulo the magnitude of a straight angle, denoted by $\varpi$. Its result is called `sum of angular magnitudes'.

	\textbf{Note 2}\quad The meaning of the expression $[\angle\,pr]_{\cong}\oplus[\angle\,rq]_{\cong}=[\angle\,pq]_{\cong}$ is the following: Take an angle $\angle\,pr$ of equivalence class $[\angle\,pr]_{\cong}$ and an angle $\angle\,rq$ of equivalence class $[\angle\,rq]_{\cong}$, both lying in the same plane and having a common vertex $O$, and join them at their common side $r$, such that $r$ lies between the other two sides $p$ and $q$, respectively. The result is an angle $\angle\,pq$, which is in the equivalence class $[\angle\,pq]_{\cong}$.

	\textbf{Note 3}\quad The expression $[\angle\,pr]_{\cong}\oplus[\angle\,rq]_{\cong}=[\angle\,pq]_{\cong}$ may be written in the form $\alpha\oplus\beta=\gamma$. However, it should be kept in mind, that the Greek letters do not denote angles, but angular magnitudes.
\end{quote}

The `addition of angular magnitudes' is often called just `addition of angles'. This should be avoided, since the two terms have a different meaning.

\section{Angular measures}

Angular magnitudes are not numbers, but mathematical objects of a different kind. However, the binary operation $\oplus$ introduced in Theorem \ref{Def3} shares some similarities with the addition of real numbers. Therefore, in order to assign real numbers to angular magnitudes, there is a strong temptation to look for a morphism between the semigroup $(\mathcal{A},\oplus)$ and the additive group of the positive real numbers $(\mathbb{R}^+,+)$. But this would not be in accordance with the fact that it is not magnitudes but rather the ratio of magnitudes of the same kind that reflects physical reality.

When dealing with magnitudes, the physically correct way to proceed is to establish a proportion\footnote{A proportion establishes an equivalence relation between two ratios. A ratio is a relation and should not be confused with a quotient, which is a binary operation defined for numbers.} between a ratio of two magnitudes of the same kind and a ratio of two numbers, {i.\,e.}\ instead of assigning real numbers directly to magnitudes, we rather assign them to their ratios.

It is known since ancient times that there is a relationship between an angle and an arc of a circle whose centre is located at the vertex of the angle subtending this arc. The Babylonians used this knowledge in a pragmatic way to define an angular measure by a uniform subdivision of the circumference of the circle. This measure is still in use today. But it was the ancient Greeks who first proved the relationship between circular arcs and angles.

In equal circles, according to Euclid \cite{Euclid1968} (\emph{Elements}, Book III, Prop.~26 and 27), angles whose vertices are subtended by equal arcs are equal and vice versa. Furthermore (\emph{Elements}, Book VI, Prop.~33\footnote{There is an inaccuracy in the proof of this theorem, since the sum of the angles can exceed two right angles. This contradicts the restriction to angles less than or equal to two right angles as given by Euclid. However, the proof can be amended (see e.g. Legendre \cite{Legendre1828} and Crelle \cite{Crelle1826}).}), the ratio of the angular magnitudes is equal to the ratio of the lengths of the arcs subtended by the vertices of the respective angles, {i.\,e.}\ the proportion
\begin{equation}
	\alpha:\beta::s_{\alpha}:s_{\beta}
\label{prop1}
\end{equation}
is valid, where $\alpha$ and $\beta$ denote angular magnitudes and $s_{\alpha}$ and $s_{\beta}$ refer to the corresponding arc lengths. In addition Archimedes of Syracuse has shown (\emph{Measurement of a Circle}, Prop.~3), that for any circle the ratio of its circumference to its diameter is a constant number, which we today denote by $\PI$.\footnote{In fact Archimedes only gave an upper and a lower limit for $\PI$, but his algorithm offers the possibility to approximate this number with any desired accuracy.} Hence we have another proportion
\begin{equation}
	c:r::2\PI:1\,,
\label{prop2}
\end{equation}
where $c$ denotes the circumference of the circle, $r$ its radius and the real number $\PI$ is a metric constant of the circle in Euclidean geometry.

If we now set in proportion \eqref{prop1} $s_{\beta}=c$ and $\beta=2\,\varpi$, {i.\,e.}\ the central angle of a full circle, and use proportion \eqref{prop2} we obtain the proportion
\begin{equation}
	\PI\,(\alpha:\varpi)::s_{\alpha}:r\,.
\label{eqn3}
\end{equation}
The multiplication of the ratio on left side of this proportion by $\PI$ means, that the real number assigned to the ratio of angular magnitudes is $\PI$ times the real number assigned to the ratio of lengths on the right side of the proportion.

Since both the left and the right side of proportion \eqref{eqn3} are ratios of magnitudes of the same kind, we can equate each side with the same real number. This yields the two equations
\begin{equation}
	s_{\alpha}=\varphi_{\alpha}\,r
\label{eqn4a}
\end{equation}
and
\begin{equation}
	\alpha=\frac{\varphi_{\alpha}}{\PI}\,\varpi\,,
\label{eqn4b}
\end{equation}
where $\varphi_{\alpha}$ denotes a positive real number, called `angular measure'. Note that $\varphi_{\alpha}$ depends neither on an ``angular unit'' nor on the choice of a length unit, because the angular measure is of dimension number (for details see \cite{Krystek2015}).

It is easy to verify that the angular measure takes values in the interval $]0,2\PI]$, where the value $2\PI$ corresponds to the central angle of a full circle.

\section{Reference angles and the circle graduation}

In practical measurement of angles, a circle graduation is used instead of the angular measure. The main reason is, that a sectioning of the angular measure leads to irrational values, which cannot be subject to measurement.

If we set in equation \eqref{eqn4b} ${2\varpi=p\,\epsilon}$, where $\epsilon$ denotes a certain fraction of the central angle of a full circle, denoted by $2\varpi$, we obtain
\begin{equation}
	\alpha=\varphi_{\alpha}\frac{p}{2\PI}\,\epsilon\,.
\label{eqn4c}
\end{equation}
In angle measurement technology, the angular magnitude $\epsilon$ is called a `reference angle'. The reference angle is often imprecisely called angular unit, but in contrast to length measurement, no unit is required for angle measurement, as was clearly shown in the previous section. Angle measurement can be accomplished without calibration. It is based on the principle of circle closure, a natural conservation law for plane angles, known since the time of Euclid, expressing the fact that the sum of the angular measures around any point in a plane is equal to $2\PI$.

If we now express the angular magnitude $\alpha$ by its numerical value $\varphi_{\alpha}$ with respect to a reference angle, we obtain
\begin{equation}
	\varphi_{\alpha}=\frac{2\PI}{p}\{\alpha\}\,,
\label{eqn5}
\end{equation}
where both $\{\alpha\}$ and $p$ denote the numerical value of an angular magnitude with respect to an arbitrary reference angle. The special value $p$ is assigned to the central angle of a full circle $2\varpi$. Note that the numerical values $\{\alpha\}$ and $p$ change, of course, with a change of the reference angle, while their quotient always remains constant, independent of the selected reference angle.

Equation \eqref{eqn4a} relates the length of a circular arc $s_{\alpha}$, subtended by the vertex of an angle with an angular magnitude $\alpha$, to the radius $r$ of that arc. The angular measure $\varphi_{\alpha}$ occurring in this equation is a constant of proportionality whose value can be calculated by using equation \eqref{eqn5}. This equation depends on the choice of a reference angle, but it is valid for any reference angle. The reference angle is implicitly fixed by the choice of the real number $p$. Note that this factor is also a constant of proportionality which in this case relates the angular measure $\varphi_{\alpha}$ to the numerical value $\{\alpha\}$ assigned to an angular magnitude $\alpha$. It has to be emphasized that both $\varphi_{\alpha}$ and $2\PI/p$ are pure real numbers.

Equation \eqref{eqn5} shows that the angular measure $\varphi_{\alpha}$ and the numerical value $\{\alpha\}$ associated with an angular magnitude $\alpha$ are generally not identical, but only proportional. However, it would be more natural, if these two values were identical. This can be achieved, if we in particular choose ${p=2\PI}$, whereby equation \eqref{eqn5} is simplified to
\begin{equation}
	\varphi_{\alpha}=\{\alpha\}\,.
\label{eqn6}
\end{equation}
In this case the numerical value $\{\alpha\}$ is related to the reference angle `radian'. This can easily be seen by combining equations \eqref{eqn6} and \eqref{eqn4a} and subsequently setting ${\{\alpha\}=1}$, which yields ${s_{\alpha}=r}$, {i.\,e.}\ the arc is measured by its radius, as it is required by the definition of the radian.

In everyday life and in technical applications, the reference angle 'degree', with symbol $^\circ$, is usually used, which corresponds to the value $p=360$.

\section{Trigonometric functions}

Let an angle with magnitude $\alpha$ be given. If a circle with radius $r$ is drawn about the vertex of this angle, it is intersected by the sides of this angle in two points. We denote the length of the arc between these two points as before by $s_{\alpha}$ and the length of the chord subtending this arc by $\cord(\alpha)$\footnote{This was the trigonometric function of the ancient Greek, which is best known from the table of chords in Book I, chapter 11, of the μαθηματική σύνταξις (also known as Almagest) written by Claudius Ptolemy during the 2$^\text{nd}$ century.}. Using modern calculus, we find that the length of the arc is related to the length of the chord by
\begin{equation}
\frac{s_{\alpha}}{2r}=F\left(\frac{\cord(\alpha)}{2r}\right)\,,
\label{eqn7}
\end{equation}
with the function
\begin{equation}
	F:[0,1]\to\mathbb{R}\,,
	\qquad
	x\mapsto\int\limits_{0}^{x}\!\!\frac{\D t}{\sqrt{1-t^2}}\,.
\label{eqn8}
\end{equation}
The domain of $F$ needs to be restricted to the interval $[0,1]$, because a chord can never become greater than the diameter of the same circle.

$F$ is monotonically increasing and maps the interval $[0,1]$ to the interval $[0,\PI/2]$. Thus there exists an inverse function which we will denote by $S$. Combining equations \eqref{eqn5} and \eqref{eqn7} and applying the inverse function $S$ yields
\begin{equation}
	\frac{\cord(\alpha)}{2r}=
	S\left(\frac{2\PI}{p}\cdot\frac{\{\alpha\}}{2}\right)\,.
\label{eqn9}
\end{equation}

If the considered angle is bisected, two similar right-angled triangles are obtained in which the side opposite to the vertex of the angle is equal to half of the chord and whose hypotenuses are equal to the radius of the circle. In trigonometry the ratio of these two lengths is called the sine of the related angle. Thus the function $S$ is a certain version of the sine function restricted to the domain $[0,\PI/2]$. By real analytic continuation of $S$ we can define the function
\begin{equation}
	\Sin_p:\mathbb{R}\to\mathbb{R}\,,
	\qquad
	\{\alpha\}\mapsto\sin\left(\frac{2\PI}{p}\{\alpha\}\right)\,.
\label{eqn10}
\end{equation}
Furthermore, by using the Pythagorean theorem we can also introduce the function
\begin{equation}
	\Cos_p:\mathbb{R}\to\mathbb{R}\,,
	\qquad
	\{\alpha\}\mapsto \cos\left(\frac{2\PI}{p}\{\alpha\}\right)\,,
\label{eqn11}
\end{equation}
such that
\begin{equation}
	(\Cos_p x)^2+(\Sin_p x)^2=1
\label{eqn12}
\end{equation}
is valid. The functions $\Sin_p$ and $\Cos_p$ can be considered as a generalisation of the sine and cosine function, respectively. They are continuous on $\mathbb{R}$ and all known algebraic rules for trigonometric functions apply. A generalisation of all other trigonometric functions can be derived from them in the usual way.

The function terms $\Sin_p$ and $\Cos_p$ are related by Euler's identity \cite{Euler1748}
\begin{equation}
	\exp\left(\I\frac{2\PI}{p}\{\alpha\}\right)=
	\cos\left(\frac{2\PI}{p}\{\alpha\}\right)+
	\I\sin\left(\frac{2\PI}{p}\{\alpha\}\right)\,.
\label{eqn12a}
\end{equation}
The function $\exp$ maps the real numbers $\mathbb{R}$ to the unit circle $\mathbb{S}^1$ in the complex plane. It is therefore relevant for the definition of the so-called `phase angle'.

We can also introduce the functions
\begin{equation}
	\Arcsin_p:[-1,1]\to\mathbb{R}\,,
	\qquad
	x\mapsto\frac{p}{2\PI}\arcsin x
\label{eqn10a}
\end{equation}
and
\begin{equation}
	\Arccos_p:[-1,1]\to\mathbb{R}\,,
	\qquad
	x\mapsto\frac{p}{2\PI}\arccos x\,,
\label{eqn11a}
\end{equation}
which are generalisations of the inverse functions of $\Sin_p$ and $\Cos_p$, respectively. Its result is a real number of the interval $[0,p]$, which is a numerical value assigned to an angular magnitude. This numerical value is called `principal measure' of an angle relative to $p$. It allows a classification of angles as given in Table \ref{tab2}.

\begin{table}[htb]
	\centering
	\caption{Classification of angles by their principal measures}
	\vspace*{0.25\baselineskip}
	\begin{tabular}{lc}
		\toprule
		\textbf{name of the angle} & \textbf{principal measure} \\
		\midrule
		zero angle & $0$ \\[1mm]
		acute angle & $\Big]0,\dfrac{p}{4}\Big[$ \\[3mm]
		right angle & $\dfrac{p}{4}$ \\[3mm]
		obtuse angle & $\Big]\dfrac{p}{4},\dfrac{p}{2}\Big[$ \\[3mm]
		straight angle & $\dfrac{p}{2}$ \\[3mm]
		reflex angle & $\Big]\dfrac{p}{2},p\Big[$ \\[3mm]
		perigon & $p$ \\[1mm]
		\bottomrule
	\end{tabular}
\label{tab2}
\end{table}

The functions $\Sin_p$ and $\Cos_p$ are periodic with the `principal period' $p$. If in particular the principal period ${p=2\PI}$ is chosen, {i.\,e.}\ if the radian is used as reference angle, we obtain by using equation \eqref{eqn6}
\begin{equation*}
	\Sin_{2\PI}\alpha=\sin\varphi_{\alpha}
	\qquad\text{and}\qquad
	\Cos_{2\PI}\alpha=\cos\varphi_{\alpha}\,.
\end{equation*}
These are the sine and cosine functions commonly used in trigonometry. The corresponding inverse functions are
\begin{equation*}
	\Arcsin_{2\PI}\alpha=\arcsin\varphi_{\alpha}
	\qquad\text{and}\qquad
	\Arccos_{2\PI}\alpha=\arccos\varphi_{\alpha}\,.
\end{equation*}
An addition we obtain for $p=2\PI$ from equations \eqref{eqn12a} and \eqref{eqn6} the relation
\begin{equation*}
	\E^{\I\varphi_{\alpha}}=\cos\varphi_{\alpha}+\I\sin\varphi_{\alpha}\,,
\end{equation*}
which is important for the `phase angle'. The phase angle used in mathematics and theoretical physics is equal to the angular measure $\varphi_{\alpha}$.

It turns out that the argument of \emph{all} commonly used trigonometric functions, as with all other transcendental functions, is \emph{not} an angle, but rather an angular measure, {i.\,e.}\ a pure number. It is therefore neither necessary nor correct to add `rad' to their argument, since this particular reference angle is already implicitly part of the definition of these functions. The results of the commonly used inverse trigonometric functions are real numbers, which are, of course, related to the reference angle radian, which is implied by the definition of these functions as well.

\section{Conclusion}

It has been shown that most misunderstandings regarding the `radian' can be avoided if some simple rules are obeyed. Although the discussion has been confined to the plane angle of Euclidean geometry, all conclusions apply equally well to the concepts of `angle of rotation' and `phase angle', which have not been discussed here, in order to concentrate on the essential points.

When dealing with the radian,
\begin{equation*}
	\hspace*{-15mm}\text{the \emph{wrong} equation}\quad
	\alpha=\frac{s}{r}
	\quad\bigl(\text{instead of the \emph{correct} equation}\quad
	s=\varphi_{\alpha}r\bigr)
\end{equation*}
is usually used in the literature. From this wrong equation then the conclusion is drawn that the radian is a ``derived unit'' which is equal to the number one. Unfortunately, this statement also appears in the current SI brochure, where moreover `rad' is expressed by the quotient m$/$m, in order to emphasize that it \emph{is} ``derived''. But these statements are both not justified at all.

If any value associated with an angular magnitude is reported, \emph{both} the numerical value and the corresponding symbol of the reference angle must \emph{always} be stated, because the numerical value depends on the chosen reference angle. In case of a semicircle, for example, the value associated with the angular magnitude has to be reported by $\PI\,\text{rad}$, \emph{not} simply by $\PI$, i.e. the symbol `rad' \emph{must not} be omitted. On the other hand, ${c=\PI\,r}$ \emph{must always} be written for the arc of a semicircle with radius $r$, i.e. in this case it is \emph{mandatory} to omit the symbol `rad', because the angular measure has to be used here, which is a pure number.

When dealing with trigonometric functions, \emph{wrong} expressions, such as $\sin\alpha$ instead of the \emph{correct} expression $\sin\varphi_{\alpha}$ are usually used in the literature. From this incorrect expressions then the conclusion is drawn that the symbol `rad' has to be appended to the numerical value of their argument. But this is wrong. The argument of the trigonometric functions is \emph{always} an angular measure, {i.\,e.}\ it is neither necessary nor correct to add the symbol `rad'.

\printbibliography

\end{document}